\title{Using Simple Linear Models with Truncation to Determine the Gregorian Day of the Week}
\author{Bryce Iversen\thanks{Sonoma State University and University of Illinois, Urbana-Champaign}}
\date{August 24, 2025}

\documentclass{article}

\usepackage{graphicx}
\usepackage{amsmath}
\usepackage{amssymb}

\newenvironment{AMS}{\em 2020 Mathematics Subject Classification:}{}
\newenvironment{keywords}{\em Keywords and phrases:}{}

\newcommand{\floor}[1]{\left \lfloor {#1} \right \rfloor}

\begin{document}

\sloppy

\newpage

\maketitle

\begin{abstract}
	The Gregorian calendar---first established for daily use on Friday, October 15th, 1582 by Pope Gregory XIII in Catholic countries---is presently the most pervasive calendar in the world. As such, algorithms for performing various calendrical computations in accurate, performant, and easily implementable ways are extremely useful in fields like software engineering. In this paper, we present a novel algorithm for determining the day of the week for any date in the Gregorian calendar. Of note, our algorithm does not rely on remembering tables of values. Instead, we encode tables needed for computation using simple linear regression with truncation to adjust for any errors present in our linear models in such a way that no tables have to be recalled. In addition, our algorithm does not require a relabeling of days, weeks, months, or years to values other than their intuitive representations. The algorithm works by taking a date in the Gregorian calendar, calculating the number of days (accounting for leap years using simple linear regression with truncation) that have elapsed since the epoch of the Gregorian calendar in 1582 from the specified date and adding this number modulo 7 to the epoch's day of the week thus, obtaining the day of the week for the requested date in numeric form.
\end{abstract}



\section{Introduction}

Calculating the day of the week in the Gregorian calendar is fraught with error. A year has 365 days---except during a leap year when it has 366. A month can have anywhere between 28 and 31 days with February having 28 days except on leap years when it has 29 days instead. A leap year occurs on years divisible by four--except for years divisible by 100 (centurial years) but those centurial years can become centurial leap years if they are also divisible by 400. Why is the Gregorian calendar such a mess? Put simply, these leap year adjustments ensure the calendar year does not fall out of sync with the solar year as centuries pass. Unfortunately, leap years and the varying length of months makes a mess out of calendrical calculations. What makes predicting the day of the week tricky is that there are 5 levels of temporal granularity that all must be properly addressed: centuries (for accurately calculating leap years), years, months, weeks, and days. Thankfully, weeks are always 7 days long and there are always exactly 12 months in a year.

To ease the confusion of the Gregorian calendar, we adopt several conventions. Calendars are to be described according to ISO 8601~\cite{iso8601}. This means, we assume weeks start on Monday and end on Sunday. More specifically, we assign the following numerical values to the days of the week to denote their order and make them usable in computations: $\text{Monday} = 0$, $\text{Tuesday} = 1$, $\text{Wednesday} = 2$, $\text{Thursday} = 3$, $\text{Friday} = 4$, $\text{Saturday} = 5$, and $\text{Sunday} = 6$.\footnote{We actually stray from ISO8601 here because ISO8601 assigns Monday to the number one and Sunday to the number seven. However, we choose the convention outlined in this paper instead for the purposes of modulo arithmetic.} Similarly, we assign the months of the year numerical values for the same reasons as above: $\text{January} = 1$, $\text{February} = 2$, $\text{March} = 3$, $\text{April} = 4$, $\text{May} = 5$, $\text{June} = 6$, $\text{July} = 7$, $\text{August} = 8$, $\text{September} = 9$, $\text{October} = 10$, $\text{November} = 11$, $\text{December} = 12$. The epoch, or the beginning of the Gregorian calendar, was Friday, October 15th, 1582. Where henceforth, this date may be encoded unambigiously as 1582-10-15---the four-digit, left-padded year separated by a hyphen followed by the two-digit, left-padded month broken by another hyphen and then trailed by the two-digit, left-padded day.

\section{Methods}
The algorithm described in this paper was originally discovered, implemented, and tested in the {\em C} programming language. All dates from 1582-10-15 to 9999-12-31 were computed with the algorithm and compared against the output of Zeller's congruence~\cite{zeller}. The algorithm was eventually ported to {\em Python} for greater portability, accessibility, and legibility. The CSV files found in the {\em GitLab} repository~\cite{repo} referenced in this paper were generated using the aforementioned {\em C} program in conjunction with the {\em LibreOffice Calc} spreadsheet software and the {\em R} programming language. All figures found in this paper were generated in the {\em R} programming language.
Throughout this paper, when referencing ``simple linear regression", we are regressing using the Oridinary Least-Squares (OLS) algorithm~\cite{anton-rorres}. All closed-forms for linear equations originating from OLS were found using linear algebra formulations on a {\em TI-Nspire\textsuperscript{TM} CAS CX II}.

\section{The Formula}

The formula we will begin to describe can be summarized as merely adding the days since the Gregorian epoch's day of the week and taking this sum modulo 7.

\begin{equation}
	W(y, m, d) = \Bigl[ E_w + 365 \left( y - E_y \right) - 1 + L(y) + D(y, m, d) \Bigr] \pmod{7}
\end{equation}

\noindent
where:

\begin{description}
	\item $W(y, m, d)$ is a function that takes in the date and computes the day of the week,
	\item $E_w$ is the Gregorian epoch's day of the week (Friday or $4$),
	\item $E_y$ is the Gregorian epoch's year ($1582$),
	\item $L(y)$ is a function that counts the number of leap-years that have occured since the Gregorian epoch,
	\item $D(y, m, d)$ is a function that counts the days into the year that the specified date occurs.
\end{description}

First, a date is encoded into a year $y$, a month $m$ and day $d$. The $365 \left( y - E_y \right) - 1$ term approximates the number of days that have occured since the Gregorian epoch up until the last day of the year $y - 1$ by calculating the number of years that have elapsed and multiplying by the number of days in a non-leap year: 365 days.

However, this term is an underestimate in general because in a leap year there are 366 days. So, we must add a day for each leap-year that has occured. This is exactly what the function $L(y)$ does. Therefore, $365 \left( y - E_y \right) - 1 + L(y)$ gives us the number of days that have elapsed from the end of the year $1582$ to the end of the year $y - 1$. We also have to count the days that have elapsed into the year $y$. The function $D(y, m, d)$ does exactly this. Hence, $365 \left( y - E_y \right) - 1 + L(y) + D(y, m, d)$ gives us the number of days that have occured since the end of $1582$ to the date specified by $y$, $m$, and $d$. Finally, adding the day of the week at the Gregorian epoch modulo 7 to this sum should give us the day of the week of the date $y$, $m$, and $d$. A keen observer might note however, that we seem to have omitted the days between the date of the Gregorian epoch to the end of $1582$. As it turns out, we don't need to count these days since 1582-10-15 occurs on a Friday and 1582-12-31 also occurs on a Friday meaning that moduolo 7, we can decide to merely add $E_w$ to $365 \left( y - E_y \right) - 1 + L(y) + D(y, m, d)$ and we will obtain the same value as if we had counted every day inbetween 1582-10-15 and 1582-12-31.

This is also part of the reason why we must subtract $1$ from $365 \left( y - E_y \right)$. Consider again Saturday, January 1, 1583. Hence, if $y = 1583$, then $365 \left( y - E_y \right) = 365 \left( 1583 - 1582 \right) = 365$. An astute observer might note however, that 365 days have not elapsed since the epoch on the day of 1583-01-01! However, we can play a modulo arithmetic trick here while preserving the function this term performs. We start by noting that $365 \left( y - E_y \right) = 365 = 365 \pmod{7} = 1 \pmod{7}$. We would like this term to be zero during the year 1583 for two reasons: (1) zero years have passed and (2) subtracting one preserves the proper day of the week for the first day of every subsequent year. To see what is happening for reason (2), observe that 1582-01-01 was a Friday, 1583-01-01 was a Saturday, and 1584-01-01 was a Sunday. As we can see, The day of the week increases by one (or sometimes two for leap years) every year on the same date. Similarly, $365 \left( 1582 - E_y \right) = 0 \pmod{7}$, $365 \left( 1583 - E_y \right) = 1 \pmod{7}$, $365 \left( 1584 - E_y \right) = 2 \pmod{7}$. Hence, $E_w + 365 \left( y - E_y \right)$ is essentially functioning as an offset for the day of the week on the first day of every year starting in 1583, when the offset is $1$ because 1582-01-01 was a Friday so 1583-01-01 should be one day later: a Saturday. The only problem with this offset is that we are also adding $D(y, m, d)$ to $E_w + 365 \left( y - E_y \right)$ which means we are double-counting the first day of the year! Since we count the first day of the year later with $D(y, m, d)$, we should subtract by one. Therefore, our final term that approximates the number of days that have passed since the year 1582 up until the end of the year $y - 1$ is given by $E_w + 365 \left( y - E_y \right) - 1$.

\subsection{Counting Leap Years}

Counting the number of leap years since a given date is tricky namely because of the exceptions to the rule of a leap year occuring every 4 years. Equation~\ref{equation:is-leap-year} is a formalization of the exact rule needed to determine if a year is a leap year. If a given year is a leap year, it will return 1 and 0 otherwise.

\begin{equation}
	\label{equation:is-leap-year}
	l(y) = 
		\begin{cases}
			1 & \text{if } \left( y \mid 4 \right) \land \left( y \nmid 100 \lor y \mid 400 \right)  \\
			0 & \text{otherwise}
		\end{cases}
\end{equation}

We consider a dataset that lists every leap year that could theoretically occur assuming the Gregorian calendar started in 1200 till the year 9999.~\cite{repo} The reason for why we will be using such a dataset is because we want to accurately predict centurial leap years using a simple linear regression model. The first centurial leap year after 1582 is 1600 but, performing OLS on a dataset consisting of the leap years since 1582 establishes an asymmetrical model fit since a centurial leap year occurs so soon after the Gregorian calendar begins. So, we extend the pattern backwards in time to 1200---the first theoretical centurial leap-year before 1582.

\begin{figure}
	\centering
	\includegraphics[width=0.5\textwidth]{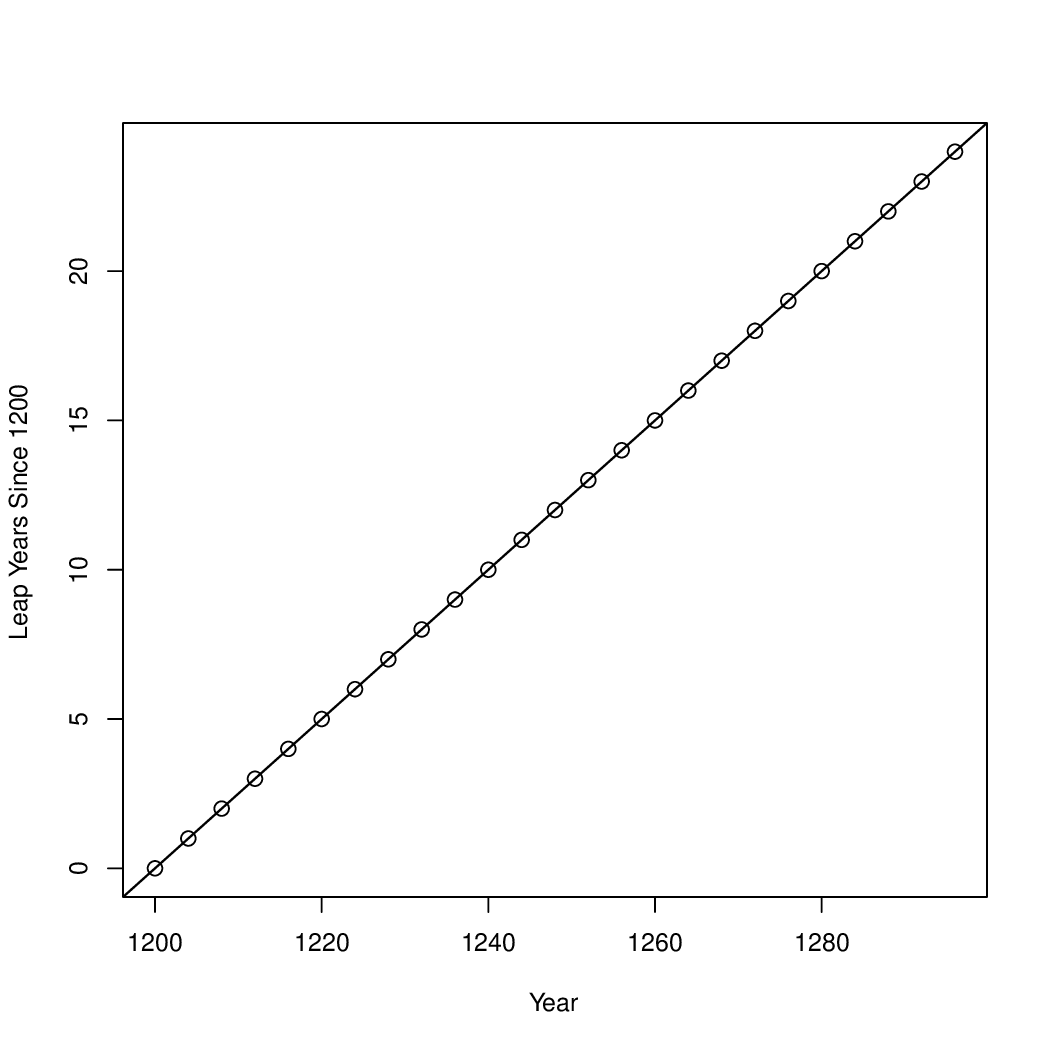}
	\caption{A plot of the leap years since 1200 up until 1300 assuming that the Gregorian method of counting leap years was extended backwards in time}
	\label{figure:1}
\end{figure}

If we observe a scatter plot of the back-dated leap years since 1200 up until 1300 (Figure~\ref{figure:1}) we can see a perfect linear trend between the year and the number of cumulative leap years ($r = 1$). The line that fits this data is given by Equation~\ref{equation:initial-leap-years-model}.

\begin{equation}
	\label{equation:initial-leap-years-model}
	\widehat{\text{LeapYearsSince1200}} = \frac{\text{Year}}{4} - 300
\end{equation}

\begin{figure}
    \centering
    \begin{minipage}{0.45\textwidth}
        \centering
        \includegraphics[width=\textwidth]{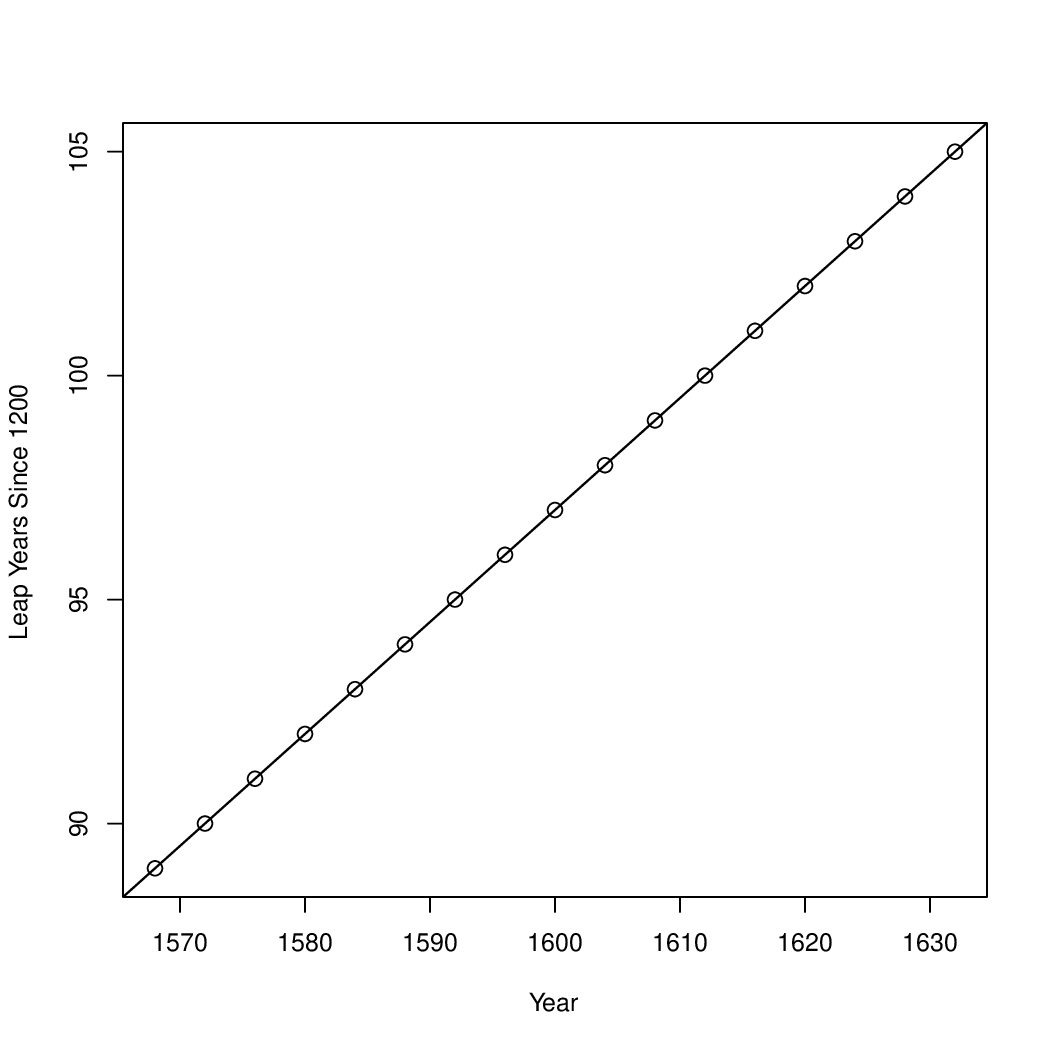}
        \caption{A plot of the leap years since 1200 from 1568 to 1632}
        \label{figure:2}
    \end{minipage}
    \hfill
    \begin{minipage}{0.45\textwidth}
        \centering
        \includegraphics[width=\textwidth]{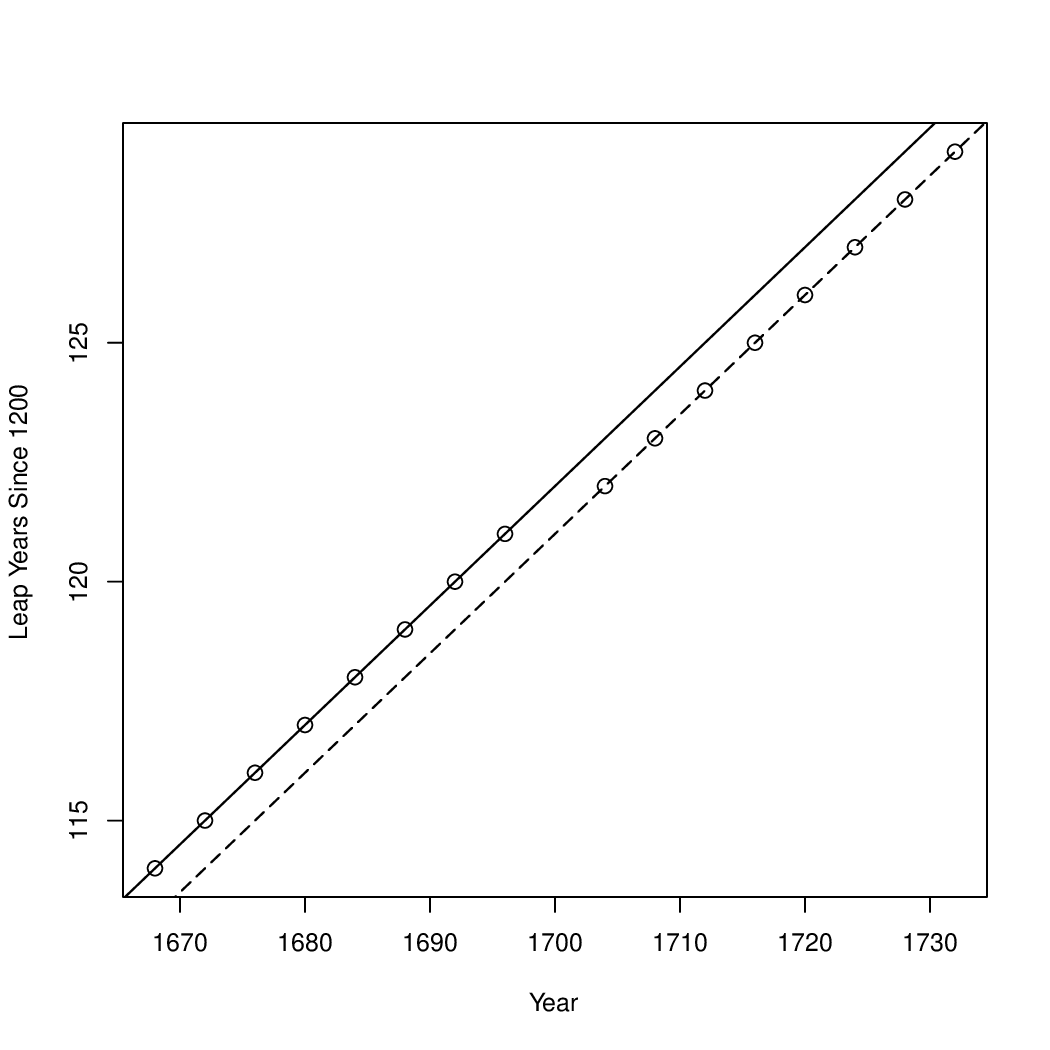}
        \caption{A plot of the leap years since 1200 from 1668 to 1732}
        \label{figure:3}
    \end{minipage}
    \vfill
\end{figure}

However, this trend is only perfectly linear if we don't look further forward in time. For example, shortly after the Gregorian calendar was established in 1582 a centurial leap year occured in 1600 (Figure~\ref{figure:2}) but did not occur a century later in 1700 (Figure~\ref{figure:3}). In Figure~\ref{figure:3} we can clearly see a gap in the data in 1700 that causes our linear model to break. But, this isn't a problem since if we compare the line fitting the leap years since 1200 from 1668 to 1696 (depicted as a solid line) to the line fitting the leap years since 1200 from 1704 to 1732 (displayed as a dashed line), we can see that the two lines have only slightly different intercepts. Additionally, they should have the same slope since leap years generally occur every 4 years. The linear equation of the solid line is given by $\widehat{\text{LeapYearsSince1200}} = \frac{\text{Year}}{4} - 303$ while the linear equation of the dashed line is given by $\widehat{\text{LeapYearsSince1200}} = \frac{\text{Year}}{4} - 304$. As expected, the two linear equations have the same slope and their intercepts only differ by $1$. This means with some minor adjustments, we can construct a simple linear model to handle these discontinuties perfectly.

\begin{figure}
	\centering
	\includegraphics[width=0.5\textwidth]{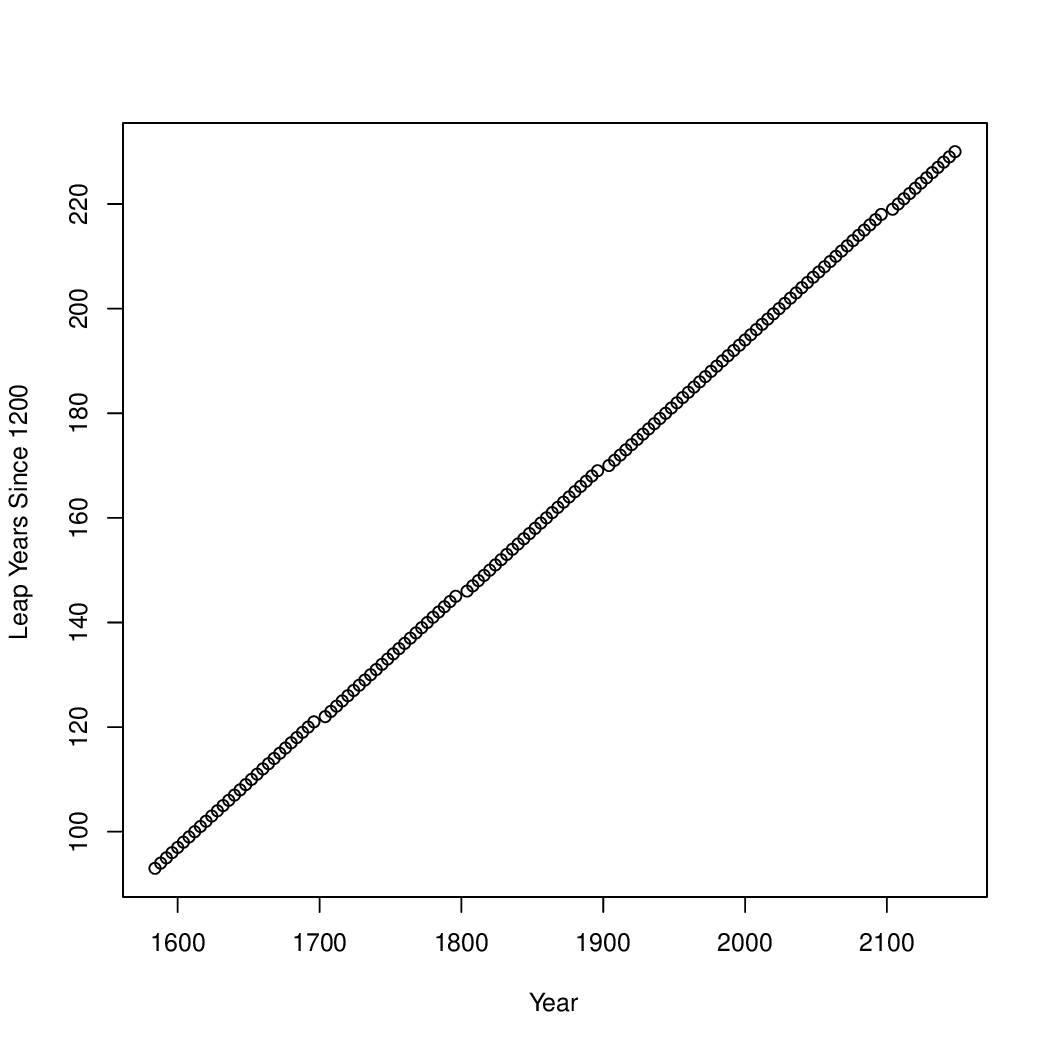}
	\caption{A plot of the leap years since 1200 from 1582 to 2150}
	\label{figure:4}
\end{figure}

Between 1200 and 1704, the intercept of an equation that perfectly fits some subset of the leap years since 1200 only changed by 3. This can be accounted for by observing that 1300, 1400, 1500, and 1700 were not centurial leap years but, 1600 was a centurial leap year. Thus, for each centurial year we must subtract 1 from the intercept and for each centurial leap year we must add one. This accurately accounts for our difference of 3. Figure~\ref{figure:4} properly demonstrates this pattern; the leap years since 1200 from 1582 to 2100 are plotted. We can see discontinuties at 1700, 1800, 1900, and 2100 but not at 1600 or 2000 since those are centurial leap years. Therefore, if we properly adjust our intercept based on the elapsed centurial years and centurial leap years, we can always fit a perfect simple linear model to the leap years since 1200. However, we are not interested in the leap years since 1200, we are interested in the leap years since 1582. By 1582, there were 92 leap years since 1200. Therefore, we just need to subtract 92 from the intercept of Equation~\ref{equation:initial-leap-years-model}, subtract $\floor{\frac{y - 1200}{100}}$ to count the centurial years since 1582, and add $\floor{\frac{y - 1200}{400}}$ to count the centurial leap years since 1582. There is one final consideration: if the input year is a leap year, we do not want to count it since for our formula it will add a day for the input year before we've counted the days into the year. In counting the days into the year, we will properly account for leap years so there is no need to do it here. Thus, we must subtract the current year if it is a leap year; subtracting $l(y)$ (Equation~\ref{equation:is-leap-year}) does just this. Equation~\ref{equation:unsimplified-leap-years-model} is the resulting unsimplified linear model which after several steps can be simplified to Equation~\ref{equation:simplified-leap-years-model}.

\begin{align}
	\label{equation:unsimplified-leap-years-model}
	L(y) &= \floor{\frac{y}{4}} - 300 - \floor{\frac{y - 1200}{100}} + \floor{\frac{y - 1200}{400}} - 92 - l(y) \\
	     &= \floor{\frac{y}{4}} - \floor{\frac{y}{100} - \frac{1200}{100}} + \floor{\frac{y}{400} - \frac{1200}{400}} - 392 - l(y) \nonumber \\
	     &= \floor{\frac{y}{4}} - \floor{\frac{y}{100} - 12} + \floor{\frac{y}{400} - 3} - 392 - l(y) \nonumber \\
	     &= \floor{\frac{y}{4}} - \floor{\frac{y}{100}} - \floor{-12} + \floor{\frac{y}{400}} + \floor{-3} - 392 - l(y) \nonumber \\
	     &= \floor{\frac{y}{4}} - \floor{\frac{y}{100}} + 12 + \floor{\frac{y}{400}} - 3 - 392 - l(y) \nonumber \\
	\label{equation:simplified-leap-years-model}
	L(y) &= \floor{\frac{y}{4}} - \floor{\frac{y}{100}} + \floor{\frac{y}{400}} - 383 - l(y)
\end{align}

\subsection{Counting the number of days into a year}

The function (Equation~\ref{equation:days-into-the-year}) for calculating the day of the year, $D(y, m, d)$, can be broken down into three parts: a linear model that predicts the day of the year (the number of days into the year) $y$ up until the first day of the month $m$, a term that corrects for the error in the aforementioned linear model, and then a final term adding the number of days into the month $m$.

\begin{equation}
	\label{equation:days-into-the-year}
	D(y, m, d) = \floor{\frac{1009m}{33} - \frac{3423}{110} + l(y)} + \Bigl( 2 - l(y) \Bigr) \floor{\frac{6}{5} - \frac{m}{10}} + d - 1
\end{equation}

\begin{figure}
	\centering
	\includegraphics[width=0.5\textwidth]{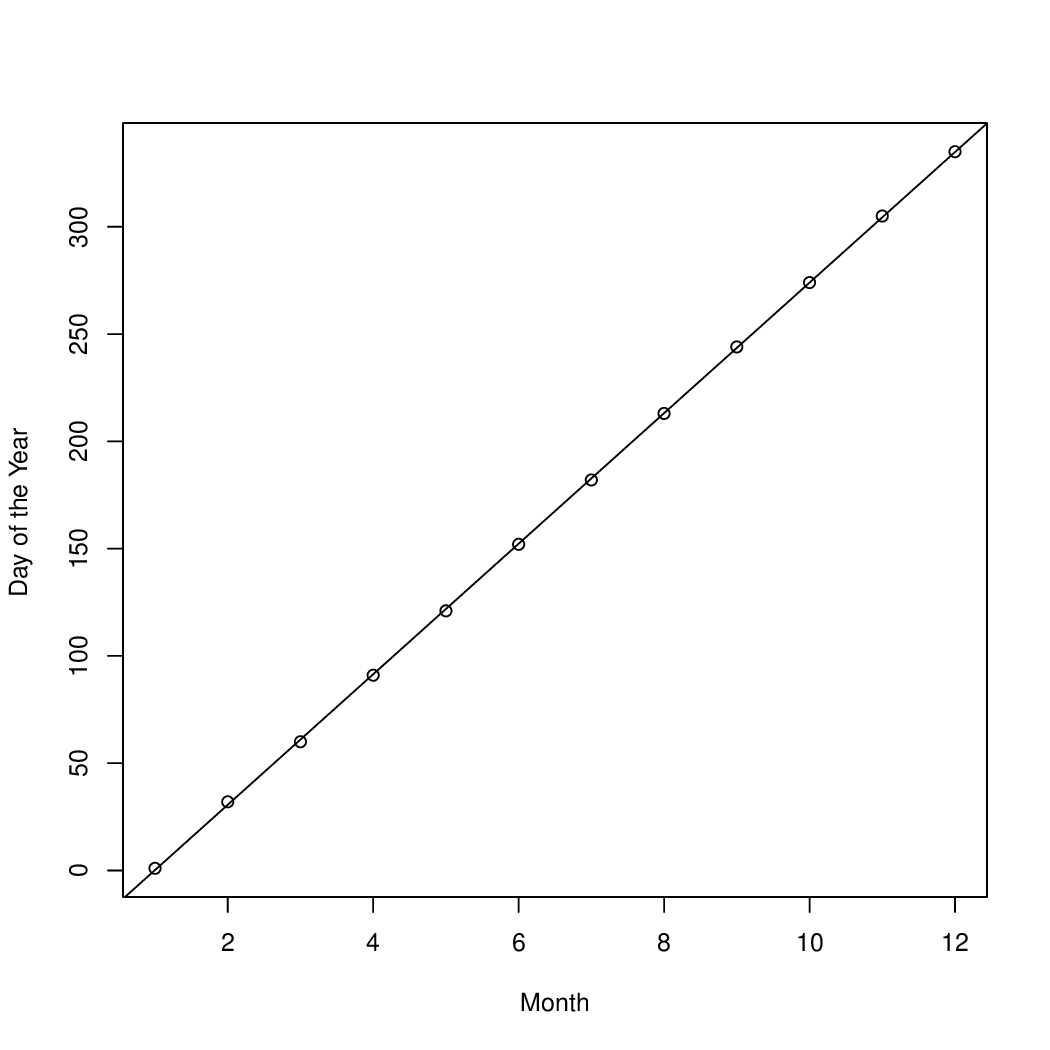}
	\caption{A plot of the linear model predicting the day of the year (for a non-leap year) from the month}
	\label{figure:5}
\end{figure}

The first term of Equation~\ref{equation:days-into-the-year}, $\floor{\frac{1009m}{33} - \frac{3423}{110} + l(y)}$, is the linear model that predicts the day of the year from the first day of a given month. This model was derived from a dataset that lists the day of the year for the first day of every month of a year that is not a leap year (the year 2001 is chosen as an example in the dataset)~\cite{repo}. Performing linear regression on this dataset with the month as the predictor and the response as the day of the year yields Equation~\ref{equation:initial-day-of-the-year-model}. Figure~\ref{figure:5} displays both the model and the linearity of the dataset.

\begin{equation}
	\label{equation:initial-day-of-the-year-model}
	\widehat{\text{Day of the Year}} = \frac{4350 \left(\text{Month} \right)}{143} - \frac{655}{22}
\end{equation}

In Table~\ref{table:initial-day-of-the-year-model-table} we can see that there is some amount of error present in the model unlike in the leap year counting model. In addition, we will need to convert the prediction to an integer at some point. Rounding, taking the ceiling, or taking the floor of the model (only rounding is shown in Table~\ref{table:initial-day-of-the-year-model-table}) will still produce some amount of error. This model is a good initial approximation, but we want a model such that when rounding to an integer by any method, it will have zero error between the prediction and the observed data.

\begin{table}
	\centering
	\caption{The initial simple linear regression model fit to the day of the year dataset (Where the year is not a leap-year)}
	\resizebox{\textwidth}{!}{
		\begin{tabular}{|l|l|l|l|l|l|}
			\hline
			Month & Day of the Year & Initial Model & Initial Model Error & Initial Model (Rounded) & Initial Model (Rounded) Error \\ \hline
			1  & 1   & 0.2   & 0.8  & 0   & 1  \\ \hline
			2  & 32  & 30.6  & 1.4  & 31  & 1  \\ \hline
			3  & 60  & 61    & -1   & 61  & -1 \\ \hline
			4  & 91  & 91.5  & -0.5 & 91  & 0  \\ \hline
			5  & 121 & 121.9 & -0.9 & 122 & -1 \\ \hline
			6  & 152 & 152.3 & -0.3 & 152 & 0  \\ \hline
			7  & 182 & 182.7 & -0.7 & 183 & -1 \\ \hline
			8  & 213 & 213.1 & -0.1 & 213 & 0  \\ \hline
			9  & 244 & 243.5 & 0.5  & 244 & 0  \\ \hline
			10 & 274 & 274   & 0    & 274 & 0  \\ \hline
			11 & 305 & 304.4 & 0.6  & 304 & 1  \\ \hline
			12 & 335 & 334.8 & 0.2  & 335 & 0  \\ \hline
		\end{tabular}
	}
	\label{table:initial-day-of-the-year-model-table}
\end{table}

Through experimentation, it can be seen that the real issue when trying to make a model fit such that we achieve nearly zero error are the months of January through March; January has 31 days, February has only 28 days (most of the time), and March has 31 days. The day of the year increases by around 30 days for each consecutive month. However, for February the day of the year only increases by 28 (and sometimes 29)---the lowest of any of the months. Trying to fit a line through these two data points surrounding February with no error is a futile goal---let alone for these points and the other nine points. Excluding the first two months of the year and performing linear regression again, we obtain a revised linear model (Equation~\ref{equation:revised-day-of-the-year-model}) that when rounded to the nearest integer, produces no error except for the months of January and February (Table~\ref{table:revised-day-of-the-year-model-table}).

\begin{equation}
	\label{equation:revised-day-of-the-year-model}
	\widehat{\text{Day of the Year}} = \frac{1009 \left(\text{Month} \right)}{33} - \frac{1739}{55}
\end{equation}

\begin{table}
	\centering
	\caption{The revised simple linear regression model fit to the day of the year dataset (Where the year is not a leap-year)}
	\resizebox{\textwidth}{!}{
		\begin{tabular}{|l|l|l|l|l|l|}
			\hline
			Month & Day of the Year & Revised Model & Revised Model Error & Revised Model (Rounded) & Revised Model (Rounded) Error \\ \hline
			1  & 1   & -1    & 2    & -1  & 2 \\ \hline
			2  & 32  & 29.5  & 2.5  & 30  & 2 \\ \hline
			3  & 60  & 60.1  & -0.1 & 60  & 0 \\ \hline
			4  & 91  & 90.7  & 0.3  & 91  & 0 \\ \hline
			5  & 121 & 121.3 & -0.3 & 121 & 0 \\ \hline
			6  & 152 & 151.8 & 0.2  & 152 & 0 \\ \hline
			7  & 182 & 182.4 & -0.4 & 182 & 0 \\ \hline
			8  & 213 & 213   & 0    & 213 & 0 \\ \hline
			9  & 244 & 243.6 & 0.4  & 244 & 0 \\ \hline
			10 & 274 & 274.1 & -0.1 & 274 & 0 \\ \hline
			11 & 305 & 304.7 & 0.3  & 305 & 0 \\ \hline
			12 & 335 & 335.3 & -0.3 & 335 & 0 \\ \hline
		\end{tabular}
	}
	\label{table:revised-day-of-the-year-model-table}
\end{table}

However, it is important to note that Equation~\ref{equation:revised-day-of-the-year-model} only predicts the day of the year for a year that is not a leap year. For a leap year, the day of the year increases by one for all months after February since February now contains 29 days instead of 28. Since the dataset we regressed on remains functionally identical and the linear model explicitly excludes January and February, all we would have to do to handle the leap year case is add one to our model when the current year is a leap year. The function $l(y)$ does exactly this. Thus, given a month $m$ and a year $y$, $\frac{1009m}{33} - \frac{1739}{55} + l(y)$ gives us the days into the year $y$ for the first of the month $m$ (with some amount of error for January and February and because we have not rounded the linear model yet). Rounding the revised model yields no error for 10 out of the 12 months, so the final modification to our linear model will be to round it. We can express the rounding function as $\text{Round}\left( x \right) = \floor{x + \frac{1}{2}}$ for some $x \in \mathbb{R}$. Thus,

\begin{align}
	\text{Round}\left( \frac{1009m}{33} - \frac{1739}{55} + l(y) \right) &= \floor{ \frac{1009m}{33} - \frac{1739}{55} +l(y) + \frac{1}{2} } \nonumber \\
	\label{equation:day-of-the-year-model-term}
	&= \floor{ \frac{1009m}{33} - \frac{3423}{110} + l(y) }\text{.}
\end{align}

Therefore, we have derived the day of the year linear model term but, as we have seen, it has some error for the months of January and February. During a non-leap year, this model predicts 2 days less than what we should observe for both January and February but only 1 day less than observed in a leap year. To correct for the error in the model, we need to find a function $f: \mathbb{R} \longrightarrow \mathbb{R}$ such that $f(1) = f(2) = 1$ and $f(3) = f(4) = ... = f(11) = f(12) = 0$, or, in other words, we need to find a function $f(m)$ such that when we substitute for a month $m$, we obtain 1 for January and February, but zero for all other months. Then, all we would have to do is multiply $f$ by $2$ or $1$ depending on whether it was a leap year to accurately account for the error caused by January and February. A linear function with a vertical axis intercept slightly above 1 and a horizontal axis intercept of 12 would be a good candidate. Observe that if $h(m) = \frac{6}{5} - \frac{m}{10}$, then $h(1) = 1.1, h(2) = 1, h(3) = 0.9, ..., h(11) = 0.1, h(12) = 0$. Taking the floor of $h$ we obtain $f$ as desired since: $\floor{h(1)} = 1, \floor{h(2)} = 1, \floor{h(3)} = 0, ..., \floor{h(11)} = 0, \floor{h(12)} = 0$. Multiplying $\floor{\frac{6}{5} - \frac{m}{10}}$ by $\Bigl( 2 - l(y) \Bigr)$ means we correct for the error by 1 if it's a leap year or 2 if it's not a leap year. Thus, adding $\Bigl( 2 - l(y) \Bigr) \floor{\frac{6}{5} - \frac{m}{10}}$ to Equation~\ref{equation:day-of-the-year-model-term} corrects the error in Equation~\ref{equation:day-of-the-year-model-term}.

The last term, $d - 1$ is the simplest to justify. Since we have already counted the number of days into the year $y$ up until the first of the month $m$ with \[\floor{\frac{1009m}{33} - \frac{3423}{110} + l(y)} + \Bigl( 2 - l(y) \Bigr) \floor{\frac{6}{5} - \frac{m}{10}}\text{,}\] the only task remaining is to count the number of days into the month $m$ that the specified date falls. This is simply given by $d$ however, we must subtract one since we already counted the first day of the month with the prior terms. Hence, we are done deriving Equation~\ref{equation:days-into-the-year}.

\subsection{Determining the Day of the Week}

Putting everything together, we obtain the following formulation (Equation~\ref{equation:unsimplified-formula}):

\begin{align}
	\label{equation:unsimplified-formula}
	W(y, m, d) &= \Bigl[ E_w + 365 \left( y - E_y \right) - 1 + L(y) + D(y, m, d) \Bigr] \pmod{7} \\ \nonumber
	       E_w &= 4 \\ \nonumber
	       E_y &= 1582 \\ \nonumber
	      l(y) &= 
	      	\begin{cases}
	      		1 & \text{if } \left( y \mid 4 \right) \land \left( y \nmid 100 \lor y \mid 400 \right)  \\
	      		0 & \text{otherwise}
	      	\end{cases}\\ \nonumber
	      L(y) &= \floor{\frac{y}{4}} - \floor{\frac{y}{100}} + \floor{\frac{y}{400}} - 383 - l(y) \\ \nonumber
	D(y, m, d) &= \floor{\frac{1009m}{33} - \frac{3423}{110} + l(y)} + \Bigl( 2 - l(y) \Bigr) \floor{\frac{6}{5} - \frac{m}{10}} + d - 1
\end{align}

Which after substituting all constants and simplifying, becomes Equation~\ref{equation:simplified-formula}.

\begin{equation}
	\label{equation:simplified-formula}
	\begin{split}
	W(y, m, d) = \Biggl[ & 365y - 577811 + \\
	& \floor{\frac{y}{4}} - \floor{\frac{y}{100}} + \floor{\frac{y}{400}} - l(y) + \\
	& \floor{\frac{1009m}{33} - \frac{3423}{110} + l(y)} + \Bigl( 2 - l(y) \Bigr) \floor{\frac{6}{5} - \frac{m}{10}} + d \Biggr] \pmod{7}
	\end{split}
\end{equation}

\section{Examples}

\subsection{November 9, 1989 \normalfont{\em The Fall of the Berlin Wall}}

Thus, we have $y = 1989$, $m = 11$, and $d = 9$. So, $l(y) = l(1989) = 0$. Hence,

\begin{align*}
	L(y) &= L(1989) \\
	     &= \floor{\frac{1989}{4}} - \floor{\frac{1989}{100}} + \floor{\frac{1989}{400}} - 383 - l(1989) \\
	     &= \floor{497.25} - \floor{19.89} + \floor{4.9725} - 383 - 0 \\
	     &= 497 - 19 + 4 - 383 \\
	L(y) &= 99 \text{.}
\end{align*}

Also,

\begin{align*}
	D(y, m, d) &= D(1989, 11, 9) \\
			   &= \floor{\frac{1009 \cdot 11}{33} - \frac{3423}{110} + l(1989)} + \Bigl( 2 - l(1989) \Bigr) \floor{\frac{6}{5} - \frac{11}{10}} + 9 - 1 \\
			   &= \floor{336.\overline{3} - 31.1\overline{18} + 0} + \Bigl( 2 - 0 \Bigr) \floor{1.2 - 1.1} + 8 \\
			   &= \floor{305.2\overline{15}} + 2 \cdot \floor{0.1} + 8 \\
			   &= 305 + 0 + 8 \\
	D(y, m, d) &= 313  \text{.}
\end{align*}

Finally,

\begin{align*}
	W(y, m, d) &= W(1989, 11, 9) \\
			   &= \Bigl[ 4 + 365 \left( 1989 - 1582 \right) - 1 + L(1989) + D(1989, 11, 9) \Bigr] \pmod{7} \\
			   &= \Bigl[ 365 \left( 407 \right) + 3 + 99 + 313 \Bigr] \pmod{7} \\
			   &= \Bigl[ 148555 + 415 \Bigr] \pmod{7} \\
			   &= \Bigl[ 148970 \Bigr] \pmod{7} \\
	W(y, m, d) &= 3 \text{.}
\end{align*}

Since $W(1989, 11, 9) = 3$ and $\text{Thursday} = 3$, then November 9, 1989 falls on a Thursday.

\subsection{July 26, 2024 \normalfont{\em The Paris 2024 Summer Olympics Opening Ceremony}}

Thus, we have $y = 2024$, $m = 7$, and $d = 26$. So, $l(y) = l(2024) = 1$. Hence,

\begin{align*}
	L(y) &= L(2024) \\
	     &= \floor{\frac{2024}{4}} - \floor{\frac{2024}{100}} + \floor{\frac{2024}{400}} - 383 - l(2024) \\
	     &= \floor{506} - \floor{20.24} + \floor{5.06} - 383 - 1 \\
	     &= 506 - 20 + 5 - 384 \\
	L(y) &= 107 \text{.}
\end{align*}

Also,

\begin{align*}
	D(y, m, d) &= D(2024, 7, 26) \\
			   &= \floor{\frac{1009 \cdot 7}{33} - \frac{3423}{110} + l(2024)} + \Bigl( 2 - l(2024) \Bigr) \floor{\frac{6}{5} - \frac{7}{10}} + 26 - 1 \\
			   &= \floor{214.\overline{03} - 31.1\overline{18} + 1} + \Bigl( 2 - 1 \Bigr) \floor{1.2 - 0.7} + 25 \\
			   &= \floor{183.9\overline{12}} + 1 \cdot \floor{0.5} + 25 \\
			   &= 183 + 0 + 25 \\
	D(y, m, d) &= 208 \text{.}
\end{align*}

Finally,

\begin{align*}
	W(y, m, d) &= W(2024, 7, 26) \\
			   &= \Bigl[ 4 + 365 \left( 2024 - 1582 \right) - 1 + L(2024) + D(2024, 7, 26) \Bigr] \pmod{7} \\
			   &= \Bigl[ 365 \left( 442 \right) + 3 + 107 + 208 \Bigr] \pmod{7} \\
			   &= \Bigl[ 161330 + 318 \Bigr] \pmod{7} \\
			   &= \Bigl[ 161647 \Bigr] \pmod{7} \\
	W(y, m, d) &= 4 \text{.}
\end{align*}

Since $W(2024, 7, 26) = 4$ and $\text{Friday} = 4$, then July 26, 2024 falls on a Friday.

\section{Discussion \& Conclusion}

The algorithm in this paper was developed without referencing other pre-existing algorithms. Nevertheless, the algorithm bears some resemblance in many places to other previously discovered formulae. One of the earliest formulations was discovered by Christian Zeller in the 19th century whose congruence features several terms that adjust for leap years, an adjustment mechanism for the first two months of the year (by relabeling months), and a term that handles the number of days in a month~\cite{zeller}. The adjustment mechanism is analogous to the error correction term in Equation~\ref{equation:days-into-the-year} and counting the days into the month is similiar to the function performed by Equation~\ref{equation:days-into-the-year}. In addition, centuries are being counted in Equation~\ref{equation:simplified-leap-years-model} in a similiar fashion to Zeller's congruence. A century later, Michael Keith and Tom Craver describe a formula that utilizes virtually identical leap year-counting terms to those found in Equation~\ref{equation:simplified-leap-years-model} with yet another mechanism for dealing with February's eccentric number of days.

The formula introduced in this paper does not require tables of values to compute the day of the week and so mirrors Zeller's congruence and the algorithm in Keith et al. However, it can hardly be ignored that the previously mentioned algorithms are far less unwieldy than the formula outlined in this paper. Indeed, Xiang-Sheng Wang~\cite{wang}, John Conway~\cite{conway}, and Lewis Carroll~\cite{carroll} all provide formulations that are quite a bit easier to compute in one's head. Therefore, in light of these limitations, the algorithm that was previously derived provides the following three-fold value proposition: only one logical conditional (Equation~\ref{equation:is-leap-year}), no tables, and a somewhat interesting statistical approach to calculating the Gregorian day of the week.

\section{Future Work}

The calendar preceeding the Gregorian calendar that was used in Europe prior to 1582 was the Julian calendar. Due to the similarities between the Gregorian and Julian calendar, the algorithm described in this paper could be extended to handle Julian calendar dates as well. However, no effort has been made to do so.

\section{Acknowledgments}

Thank you to Dr. Martha Shott for her valuable feedback and effective guidance throughout the writing and publication of this paper.

\section*{About the Author}
Bryce Iversen was a Ronald E. McNair Scholar and California State University Trustees' Scholar who graduated with two baccalaureate degrees in pure and applied mathematics with a minor in philosophy at Sonoma State University. He is currently pursuing a doctorate in mathematics at the University of Illinois, Urbana-Champaign where his primary research interests are in low-dimensional topology and differential geometry. He specializes in utilizing computational techniques to visualize mathematical abstractions.

\subsection*{Bryce Iversen}
Department of Mathematics \& Statistics,
Sonoma State University,
1801 East Cotati Ave,
Rohnert Park, CA 94928

\break
\noindent
Department of Mathematics,
University of Illinois, Urbana-Champaign,
273 Altgeld Hall,
1409 West Green Street,
Urbana, IL 61801
\end{document}